\documentclass{article}
\usepackage{amsmath}

\usepackage{amssymb}

\newcommand{\CC}{\mathbb C}
\newcommand{\ZZ}{\mathbb Z}
\newcommand{\fg}{\mathfrak g}

\newcommand{\End}{\mathop{\rm End\,}\nolimits}
\newcommand{\Res}{\mathop{\rm Res\,}\nolimits}

\newenvironment{proof}[1][Proof]
           {\medbreak\noindent \emph{#1: \enspace}}
           {\quad $\square$ \par \medbreak}

{\renewcommand{\theremark}{\arabic{chapter}.\arabic{remark}}
}

\newenvironment{noindentlist}%
    {
    \begin{list}{}
        { \settowidth{\labelwidth}{}
          \leftmargin=0pt }}
    {\end{list}}

\begin{document}

\textbf{Corrections to the book ``Vertex algebras for beginners'', 
  second edition, by Victor Kac}.

\vspace{3ex}

\begin{noindentlist}{}{}
\item  p. 39, $\ell$. 3$\uparrow$; p. 49, $\ell$. 11$\uparrow$; p. 56,
$\ell$. 10$\uparrow$:~~should be $n \gg 0$ instead of $n \fg 0$

\item p. 50, $\ell$. 5:~~should be $N \gg 0$ instead of $N \fg 0$

\item p. 56, $\ell$. 3$\uparrow$ reads:~~Now, choose a system of generators 
$\{a^{\alpha}\}_{\alpha \in I} $ of $R$ viewed as a $\CC
[\partial ]-$

\item  p. 57, $\ell$. 1:~~should be $\{ a^{\alpha}_j | \alpha \in I, j \geq n 
\}$

\item  p. 67. $\ell$. 8$\uparrow$ reads:~~$= \sum^{n+1}_{i=1} (-1)^{i+1} a_{i
  \lambda_i} \gamma_{\lambda_1, \ldots, \hat{\lambda}_i, \ldots, 
\lambda_{n+1}} (a_1, \ldots, \hat{a}_i, \ldots, a_{n+1})$

\item  p. 81:~~before formula (3.1.1) a line is missing:~~Note
  that in the expansion\break(cf. (2.2.5))

\item  p. 100, $\ell$. 13:~~should be $(\varphi | \psi)$ instead of $(a|b)$

\item  p. 102, $\ell$. 2:~~should be $\varphi (w)$ instead of $\varphi (z)$

\item  p. 104, $\ell$. 5:~~$V$ should be replaced by $U$

\item  p. 118, $\ell$. 5:~~should be $(b_{(k+j)}c)$ instead of $b_{(k+j)}c$

\item  p. 118, $\ell$. 6:~~should be $(a_{(m+j)}c)$ instead of $a_{(m+j)}c$

\item  p. 191, $\ell$. 2$\uparrow$:~~should be $1- \alpha \partial^3 \nu$
instead of  $1+ \alpha \partial^3 \nu$

\item  pp. 130--131:~~Theorem~4.11 and
Proposition~4.11 are false.  A corrected version of
Section~4.11 is given below.

\end{noindentlist}

\vspace{2ex}
\begin{center}
  \textbf{4.11.~~Field algebras}

\end{center}

\vspace{3ex}

Field algebras generalize vertex algebras in the same way as
unital associative algebras generalize unital commutative
associative algebras.

A \emph{field algebra} $V$ is defined by the same data as a vertex
algebra, but weaker axioms (cf. Proposition~4.8(b)):

\begin{list}{}{}

\item (partial vacuum): $Y (|0 \rangle , z) =I_V, \, a_{(-1)} |0 \rangle
   = a$,

\item ($n$-th product): $Y (a_{(n)}b,z) = Y(a,z)_{(n)}  Y(b,z) , \, n
  \in \ZZ$.

\end{list}

Note that the $n$-th product axiom is nothing else but Borcherds
identity in the form (4.8.1) for $F=(z-w)^n$.  As in the
proof of Theorem~4.8, it follows that (4.8.1)
holds for $F=z^m$ with $m \in \ZZ_+$.  Hence the $n$-th product
axiom implies (4.6.4) for $m \in \ZZ_+$, and in
particular, the axiom (C3) of conformal algebra.

As in the case of vertex algebra, the axioms of a field
algebra imply:
\begin{displaymath}
\shoveleft{(4.11.1)  \hspace{.75in}  Y(a,z) |0 \rangle |_{z=0} = a \, , \quad Y( |0 \rangle , z) = I_V \, ,}  \hspace{1.5in} 
\end{displaymath}
 \begin{displaymath}
\shoveleft{(4.11.2)  \hspace{.75in}   Y(Ta , z) = \partial Y (a,z) = [T,Y(a,z)] \, ,}\hspace{1.5in} 
\end{displaymath}
where $T \in \End V $ is defined by $Ta = a_{(-2)} |0 \rangle$.  The
$n$-th product axiom for $n >>0$ implies \emph{weak locality}:
\begin{displaymath}
\shoveleft{(4.11.3)  \hspace{.5in}  \Res_z (z-w)^N [Y(a,z) , Y(b,w)] = 0 
      \hbox{ for } N >>0 \, .}\hspace{2in} 
\end{displaymath}
Note that weak locality of fields $a(z)$ and $b(z)$ means that
$a(z)_{(n)} b(z) =0$ for $n \geq N$, some $N$.  (Unlike the usual
locality, this is not a symmetric property.)  Then, clearly, $(za
(z))_{(n)} b(z)=0$ for $n \geq N$.  Using this remark, one can
extend the proof of Dong's lemma to the weakly local case
(assuming that all ordered pairs are weakly local).
\vspace{2ex}

\textsc{Example}~~4.11.~~Recall that any two local fields satisfy the skewsymmetry
relation (3.3.6).  This, however, fails for weakly local
fields.  In order to construct a counterexample, consider the
free bosonic field $\alpha (z) = \sum_{n \in \ZZ} \alpha_n
z^{-n-1}$ (cf. Example~3.5), and let $\beta (z) =
\sum_{n>0} n^{-1} \alpha_n z^{-n}$.  Then we have:
\begin{displaymath}
  [ \alpha (z) , \beta (w) ] = i_{w,z} (z-w)^{-1} \, .
\end{displaymath}
Hence for $j \in \ZZ_+$ we have:
\begin{displaymath}
  \alpha (z)_{(j)} \beta (z) =0 \, , \, \beta (z)_{(j)} 
  \alpha (z) =  \delta_{j0} \, .
\end{displaymath}
Therefore both pairs $(\alpha , \beta)$ and $(\beta , \alpha)$
are weakly local, but (3.3.6) fails for $a= \alpha$, $b
= \beta$, $n=0$.

  Recall that the $-1$st product axiom means:
  \begin{displaymath}
\shoveleft{(4.11.4)  \hspace{.75in}    Y(a_{(-1)} b,z )= : Y(a,z) Y(b,z): \, .} \hspace{1.5in} 
  \end{displaymath}
Replacing $a$ by $T^na$ and using (4.11.2), we see that
(4.11.4) implies the $n$-th product axiom for $n<0$.

Multiplying both sides of the $n$-th product axiom by
$(-w)^{-n-1}$ and taking summation over $n \in \ZZ$, we obtain
its equivalent form in the domain $|z|>|w|$:
\begin{eqnarray*}
\shoveleft{(4.11.5)  \hspace{.5in}   Y (Y (a,z)b,-w)c = Y (a,z-w)Y (b,-w)c}\hspace{1.75in}  \\
\nonumber   \hspace{.75in}    -  p(a,b) Y(b,-w) \sum_{j \geq 0} \partial^j_w \delta (z-w)
     \Res_x x^{(j)} Y(a,x)c \, .\hspace{1.5in} 
\end{eqnarray*}
This is immediate by the following special case of Taylor's
formula in the domain $|z|>|w|$:
\begin{displaymath}
  i_{w,x} \delta ((w+x)-z) = \sum_{j \geq 0} x^{(j)} 
  \partial^j_w \delta (z-w) \, .
\end{displaymath}
Formula~(4.11.5) implies the \emph{associativity}
property in the domain $|z| > |w|$:
\begin{displaymath}
(4.11.6)    (z-w)^N Y (Y(a,z)b,-w)c = (z-w)^N 
  Y (a,z-w)Y (b,-w)c \hbox{ for } N \gg 0 \, .
\end{displaymath}
As in Section~1.4, it is easy to show that all
holomorphic field algebras are obtained by taking a unital
associative algebra $V$ and its derivation $T$, and letting
\begin{displaymath}
  Y (a, z) b = e^{z T} (a) b, \quad a, b \in V.
\end{displaymath}

The general linear field algebra $g \ell f (U)$ defined in
Section~3.2 is not a field algebra since the field
property
\begin{displaymath}
\shoveleft{(4.11.7)  \hspace{1in}     a_{(n)} b = 0 \quad
  \mbox{for} \quad n \gg 0 }\hspace{1.75in} 
\end{displaymath}
fails in general.  However, if we take a collection of mutually weakly local fields $\{ a^\alpha
(z) \} \subset g \ell f (U)$, they generate a linear field algebra
which is a field algebra.  The $n$-th product axiom for $n \geq 0$ is
implied by (3.3.7).  Next, it is immediate to check
(4.11.1) and (4.11.2).  Weak locality is proved in
the same way as Proposition~3.2.  The $n$-th product
axiom for $n<0$ follows from (4.11.4) as explained
above.  Finally, the $-1$st product axiom is
checked by a direct calculation.

We have the following field algebra analogs of the uniqueness and 
existence theorems (obtained jointly with Bojko Bakalov).

\textsc{Theorem}~4.11.~(a)~~\emph{Let $V$ be a field algebra.  For each field
 $Y(a,z)$ define the ``opposite'' field $X(a,z)$ by the formula
 (cf. (4.2.1)):}
 \begin{displaymath}
(4.11.8) \shoveleft{  \hspace{1in}   X(a,z) b= p(a,b) e^{zT} Y (b,-z) a \, . }\hspace{1.75in}
 \end{displaymath}
\emph{Let $B(z)$ be a field which is mutually local with all fields
$X(a,z), a \in V$, on any $v \in V$, i.e.,~}
\begin{displaymath}
  (z-w)^N [B(z), X(a,z)]v =0 \hbox{ for } N \gg 0 \, .
\end{displaymath}
Suppose that (4.4.1) holds for some $b \in V$.  Then
$B(z)=Y(b,z)$.

 (b)~~ \emph{Let $V$ be a vector superspace, let $|0 \rangle$  be an
 even vector and $T$ an even endomorphism of $V$.  Let $\{
 a^{\alpha} (z) \}_{\alpha \in A}$ and $\{ b^{\beta}(z)_{\beta
   \in B}$ ($A,B$ index sets) be two collections of fields such
 that each of them satisfies conditions (i)--(v) of
 Theorem~4.5 except that in (iii) ``local'' is replaced
 by ``weakly local''.  Suppose, in addition, that all pairs
 $(a^{\alpha}(z), b^{\beta}(z))$ are local on any $v \in V$.
 Then formula (4.5.1) defines a unique structure of a
 field algebra on $V$ such that $|0 \rangle$ is the vacuum vector,
 $T$ is the infinitesimal translation operator and
 (4.5.2) holds.  The same conclusion holds if the family 
 $\{a^{\alpha}(z) \}$ is replaced by the family $\{ b^{\beta}
   (z)\}$ and the fields $Y$ are replaced by the fields $X$.
 The two field algebra structures on $V$ are related by
 (4.11.8).}

\begin{proof}
It is similar to that of Theorems~4.4 and 4.5 
  using the observation that the associativity property
  (4.11.6) is equivalent to the locality of the pair
  $(Y(a,z), X(c,z))$ on $b \in V$.
\end{proof}

Taking in this theorem all fields of a field algebra, we obtain
the following corollary.

\textsc{Corollary} 4.11.~(a)~~Vacuum and translation covariance axioms along with 
 weak locality (4.11.3) and associativity
 (4.11.6) form an equivalent system of axioms of a field 
 algebra.

(b)~~If $(V , |0 \rangle , T,Y (a,z))$ is a field algebra,
 then $(V, |0 \rangle , T,X (a,z))$, where $X(a,z)$ are defined by
 (4.11.8), is  a field algebra as well.

\textsc{Remark} 4.11a.~~It follows from the above discussion that a field algebra with
$n$-th products for $n \in \ZZ_+$ and $\partial = T$ satisfies
all axioms of a conformal algebra, except the skewsymmetry
axiom (C2), which may fail in view of Example~4.11.

\textsc{Remark} 4.11b.~~It follows from the proof of Proposition~3.3(b) that two weakly
local fields $a(z)$ and $b(z)$ for which the skewsymmetry
property (3.3.6) holds, are local.  Hence a field algebra
satisfying the skewsymmetry property (4.2.2) is a vertex algebra.
This follows also from Corollary~4.11(a).

\end{document}